\documentclass[12pt]{article}
\usepackage[utf8]{inputenc}
\usepackage{amsmath,amssymb}
\usepackage{amsthm} 
\usepackage{subfigure}
\usepackage{graphicx}
\usepackage{mathtools}
\usepackage{tikz}
\usepackage{bbm}
\usetikzlibrary{arrows}
\usetikzlibrary{intersections}
\usepackage{tkz-euclide}
\newtheorem{theorem}{Theorem}[section]

\newtheorem{lemma}[theorem]{Lemma}
\newtheorem{corollary}[theorem]{Corollary}

\usepackage{hyperref}
\begin{document}
\title{\bf A new axiom system for matroids
        \\ 1. Uniform matroid recognition}
\date{}
\maketitle
\begin{center}
\author{
{\bf Brahim Chaourar} \\
{Department of Mathematics and Statistics,
\\Imam Mohammad Ibn Saud Islamic University (IMSIU)
\\P.O. Box 90950, Riyadh 11623,  Saudi Arabia}
}
\end{center}

\begin{abstract} \noindent In this paper, we give a new axioms system based on nonseparable flats with their ranks to define a matroid. We deduce a polynomial time algorithm for deciding if a given matroid (respectively, arbitrary structure) is an uniform matroid. This problem is intractable if we use an independence or an equivalent oracle.
\end{abstract}


\section{Introduction}

We follow Oxley \cite{Oxley:1992} about matroid theory.
\\ In the next paragraphs till the definition of the mincover function, all intoduced definitions are similar to those in matroid theory but in the general case when the function $f$, which plays the role of the rank function, does not satisfy any rank function properties. So readers who are familiar with matroid theory can skip all these paragraphs.
\\ Let $E$ be a finite set and $f$ a nonnegative integer function defined on $2^E$ such that $f(\O)=0$. $(E, f)$ is called a system. We define its dual as the system $(E, f^*)$, with $f^*(X)=|X|-f(E)+f(E\backslash X)$. Note that $f(\O)=0$ ($f^*(\O)=0$, respectively) if only if $|E|=f(E)+f^*(E)$. So all these three conditions are equivalent.
\\ Let $X\subseteq E$. We say that $X\subseteq E$ is an $f$-loops if $|X|\geq 1$ and $f(X)=0$. We say that $X$ is an $f$-coloops if it is an $f^*$-loops, i.e., $|X|\geq 1$ and $f^*(X)=0$. $(E, f)$ is called $f$-simple if it is $f$-loopless, and $f(X)=1$ if and only if $|X|=1$. It is $f$-cosimple if it is $f^*$-simple. Note that $(E, f)$ is $f$-simple and $f$-cosimple if and only if $(E, f^*)$ is too. We use the notation: $\mathcal E(X)=\{ \{e\}$ such that $e\in X\}$, and $\mathcal E=\mathcal E(E)$.
\\ For a more simple presentation of the definitions and results, we consider only $f$-simple and $f$-cosimple systems $(E, f)$. It is not difficult to generalize all these in the largest general case when $(E, f)$ is not $f$-simple nor $f$-cosimple. Under our assumption that $(E, f)$ is $f$-simple and $f$-cosimple, $f(X)\geq 1$ if and only if $X\neq \O$, or equivalently, $f(X)=0$ if and only if $X=\O$, and $f(X)\geq 2$ if and only if $|X|\geq 2$, or equivalently, $f(X)=1$ if and only if $|X|=1$.
\\ We say that $X$ is an $f$-flat if $f(X)<f(Y)$ for any $Y\supset X$.  Note that, under our assumption, all subsets of $\{ \O, E\}\cup \mathcal E$ are $f$-flats. We define the collection of $f$-closures of a subset $X\subseteq E$, denoted by $cl(f, X)$, or $cl(X)$ for short, as the $f$-flats $Y\supseteq X$ such that $f(Y)=f(X)$. Note that $cl(X)$ may contain more than one subset in the general case. We use the notation $cl(X)$ for a subset or for a collection of subsets, according to the context.
\\ $X$ is $f$-separable if there exists a nonempty subset $A\subset X$ such that $f(X)=f(A)+f(X\backslash A)$. Otherwise, $X$ is $f$-nonseparable. The collection of $f$-nonseparable $f$-flats is denoted by $\mathcal L(f)$. Note that, under our assumptions, $\mathcal E\subseteq \mathcal L(f)$. Moreover, if $L\in \mathcal L(f)\backslash \mathcal E$, then $2\leq f(L)\leq |L|-1$.
\\ $X$ is $f$-locked if it is an $f$-nonseparable $f$-flat and $E\backslash X$ is an $f^*$-nonseparable $f^*$-flat. We denote by $\mathcal L(ff^*)$ the collection of $f$-locked subsets. It is clear that $\mathcal E\subseteq \mathcal L(ff^*)$ if $(E, f)$ is simple and cosimple. Moreover, if we denote by $\mathcal A^*$ the collection of subsets $E\backslash A$ such that $A\in \mathcal A$, then $\mathcal L(ff^*)=\mathcal L(f)\cap (\mathcal L^*(f^*))^*$.
\\ We say that $X$ is an $f$-independent set if $f(X)=|X|$. If it is not, then it is $f$-dependent. The collection of $f$-independent ($f$-dependent, respectively) sets is denoted by $\mathcal I(f)$ ($\mathcal D(f)$, respectively).
\\ A nonempty subset $C\subseteq E$ is called an $f$-circuit if $f(C)=|C|-1$ and $f(X)=|X|$ for any $X\subset C$, i.e., $C$ is a minimal $f$-dependent. We denote by $\mathcal C(f)$ the collection of $f$-circuits.
\\ If there is no ambiguity, we remove $f$ from all previous notations.
\\ Given a finite set $E$ and a collection of subsets $\mathcal I\subseteq 2^E$, we say that $(E, \mathcal I)$ is an independence system if the following independence axioms are satisfied:
\\ (I1) $\O\in \mathcal I$;
\\ (I2) If $I\in \mathcal I$ and $J\subseteq I$ then $J\in \mathcal I$.
\\ Given a finite set $E$ and a collection of subsets $\mathcal C\subseteq 2^E$, we say that $(E, \mathcal C)$ is a matroid whose collection of circuits is $\mathcal C$ if the following circuits axioms hold:
\\ (C1) $\O \notin \mathcal C$;
\\ (C2) If $C_1, C_2\in \mathcal C$ and $C_1\subseteq C_2$ then $C_1=C_2$;
\\ (C3) If $C_1, C_2\in \mathcal C$, $C_1\neq C_2$, and $e\in C_1\cap C_2$ then there exists $C_3\in \mathcal C$, such that $C_3\subseteq (C_1\cup C_2)\backslash \{ e\}$.
\\ Let $X$ and $Y$ two subsets of $E$. We say that $(X, Y)$ is $f$-submodular if $f(X\cap Y)+f(X\cup Y)\leq f(X)+f(Y)$. We say that $\mathcal A\subseteq 2^E$ is $f$-submodular if, for any $X, Y\in \mathcal A$, $(X, Y)$ is $f$-submodular. Moreover we say that $A\subseteq E$ is $f$-submodular if $(A, X)$ is submodular for any $X\subseteq E$.
\\
\\ Next we define what we call the mincover and minpartition functions. Let $(E, f)$ be a system, and $\O\notin \mathcal A\subseteq 2^E$ a nonempty collection of subsets of $E$. Given $X\subseteq E$, we say that $\mathcal A_X=\{ A_1, ..., A_k\}$ is an $\mathcal A$-cover ($\mathcal A$-partition, respectively) of $X$ if $X\notin \mathcal A_X$, $\mathcal A_X\subseteq \mathcal A$, and $X\subseteq \bigcup_{i=1}^k A_j$ ($X=\bigcup_{i=1}^k A_j$ and the $A_j$'s are pairwise disjoint, respectively). In both cases, $k$ is called the size of this $\mathcal A$-cover ($\mathcal A$-partition, respectively). The $f$-mincover ($f$-minpartition, respectively) function (with the respect of $\mathcal A$) is defined as follows:
\\ $f(X, \mathcal A)=Min\{ \sum\limits_{j=1}^{k} f(A_j)$ such that $\{ A_1, ..., A_k\}$ is an $\mathcal A$-cover of $X\}$;
\\ $\hat f(X, \mathcal A)=Min\{ \sum\limits_{j=1}^{k} f(A_j)$ such that $\{ A_1, ..., A_k\}$ is an $\mathcal A$-partition of $X\}$.
\\ These two functions are inspired from the max-min relation that we have between the independent set polytope of a matroid and its dual when taking the primal objective function as the characteristic vector of $X$. In the case of a matroid, both functions give the rank of $X$.
\\ We say that $(E, f, \mathcal A)$ is $f$-mincover ($f$-minpartition, respectively) bounded if $f(X)\leq f(X, \mathcal A)$ ($f(X)\leq \hat f(X, \mathcal A)$, respectively) for any $X\subseteq E$. In this case, we denote by $\mathcal A_f$ the collection of subsets $A\in \mathcal A$ for which the inequality is strict. If $\mathcal E\subseteq \mathcal A$, then both functions are cardinality bounded ($f(X, \mathcal A), \hat f(X, \mathcal A)\leq |X|$, for any $X\subseteq E$) because we can cover and partition any subset by its elements, i.e., $\mathcal E(X)$ is an $\mathcal A$-cover and also an $\mathcal A$-partition of $X$.
\\ Let $\mathcal B$ be an $\mathcal A$-cover ($\mathcal A$-partition, respectively) of $X$. We say that $\mathcal B=\{ B_1, ..., B_k\}$ is an optimal $\mathcal A$-cover ($\mathcal A$-partition, respectively) for $X$ if $f(X, \mathcal A)$ ($\hat f(X, \mathcal A)$, respectively) $=\sum_{i=1}^k f(B_i)$.
\\ Let $f_i$, $i=1, 2$, be two nonnegative integer functions. We use the notation: $f_1\leq f_2$ on $\mathcal A$ if $f_1(A)\leq f_2(A)$ for any $A\in \mathcal A$. If no collection of subsets is specified, then $f_1\leq f_2$ on $2^E$.
\\ Let $\alpha$, $n$ and $r$ be three nonnegative integers, such that $0\leq r\leq n$, $E$ a finite set on $n$ elements, $f\in 2^{E}$ a nonnegative integer function such that $r=f(E)$, and $\mathcal L\subseteq 2^{E}$. We say that $(E, f, \mathcal L)$ is a $w$-system on $n$ elements and rank $r$ if the following axioms hold:
\\ (L0) $f(\O)=0$, and $\O \notin \mathcal L$;
\\ (L1) $(E, f)$ is $f$-simple and $f$-cosimple, and $\mathcal E\subseteq \mathcal L$ (trivial nonseparable flats);
\\ (L2) For any $\O \neq X\subseteq E$, $f(X)\leq f(X, \mathcal L)$ (nonseparable flats if the inequality is strict);
\\ If, in addition, we have:
\\ (L3) $f(L_1\cap L_2)+f(L_1\cup L_2)-f(L_1)-f(L_2)\leq \alpha$ for any distinct subsets $L_1, L_2\in \mathcal L_f$ ($\alpha$-submodularity restricted to $f$-nonseparable $f$-flats);
\\ then it is called a $w_{\alpha}$-system. Moreover, a $w$-system ($w_{\alpha}$-system, respectively) $(E, f, \mathcal L)$ is called minimal if $\mathcal L=\mathcal L_f$.
\\ Now a $w$-system ($w_{\alpha}$-system, respectively) $(E, f, \mathcal L)$ is locked if, in addition,
\\ (L4) $(E, f^*, \mathcal L*)$ is again a $w$-system ($w_{\alpha}$-system, respectively).
\\ It is minimal if $\mathcal L=\mathcal L_{ff^*}$, where $\mathcal L_{ff^*}$ is the collection of nonempty subsets $X\subseteq E$ such that $X\in \mathcal L_f$ and $E\backslash X\in \mathcal L^*_{f^*}$, i.e., nonempty subsets $X\subset E$ such that $1\leq f(X)<f(X, \mathcal L)$ and $1\leq f^*(E\backslash X)<f^*(X, \mathcal L^*)$. In other words, $\mathcal L_{ff^*}=\mathcal L_f\cap (\mathcal L^*_{f^*})^*$. It is clear that $(E, f, \mathcal L)$ is a locked $w$-system ($w_{\alpha}$-system, respectively) if and only if $(E, f^*, \mathcal L^*)$ is too. Moreover, $\mathcal L_{ff^*}=(\mathcal L^*_{f^*f})^*$.
\\ We prove that (simple and cosimple) $w_0$-systems, as well as locked ones, are (simple and cosimple) matroids. Actually, general minimal (locked) $w_0$-systems are exactly general matroids.
\\ As an application, we provide several polynomial algorithms for recognizing uniform matroids. In general, matroid recognition (MR) is intractable (see \cite{RobinsonWelsh:1980}). Few studies has been done for this problem. Provan and Ball provide an algorithm for testing if a given clutter $\Omega$, defined on a finite set $E$, is the class of the bases of a matroid \cite{ProvanBall:1988}, with running time complexity $O(|\Omega|^3 |E|)$. Spinrad \cite{Spinrad:1991} improves the running time to $O(|\Omega|^2 |E|)$. In this paper, we give several polynomial algorithms on $|E|$ for recognizing uniform matroids (UMR), by varying the input structure. UMR is intractable if we use an independence or an equivalent oracle \cite{JensenKorte:1982}.
\\ The remainder of the paper is organized as follows: in Section 2, we prove that minimal $w_0$-systems are exactly matroids. Then, in Section 3, we provide three polynomial algorithms for UMR. Each one of them considers a type of structure for the input: starting from a matroid, and ending to a basic structure. Finally, we conclude in Section 4.

\section{$w_0$-systems are matroids}

First we state some properties of mincover and minpartition functions.

\begin{lemma}\label{mincover}
  Let $E$ be a finite set, $f$, $f_1$, $f_2$ three nonnegative integer functions defined on $2^E$, $\mathcal E\subseteq \mathcal A, \mathcal A_i\subseteq 2^E$, $i=1, 2$, and $\O \neq X\subseteq Y\subseteq E$. Then
\\ (i) If $\mathcal A_1\subseteq \mathcal A_2$ then $f(X, \mathcal A_2)\leq f(X, \mathcal A_1)$;
\\ (ii) If $f_1\leq f_2$ on $\mathcal A$ then $f_1(X, \mathcal A)\leq f_2(X, \mathcal A)$;
\\ (iii) $f(X, 2^E)\leq f(X, \mathcal A)\leq f(X, \mathcal A_f)$;
\\ (iv) $f(X, 2^E)\leq \hat f(X, 2^E)$;
\\ (v) $f(X, \mathcal A)\leq |X|$;
\\ (vi) $f(X, \mathcal A)\leq f(Y, \mathcal A)$.
\end{lemma}

\begin{corollary}\label{SubMinimalWsystem}
Let $(E, f, \mathcal B)$ be a minimal $w$-system, and $\mathcal E\subseteq \mathcal A\subseteq \mathcal B\subseteq 2^E$. Then $(E, f, \mathcal A)$ is a minimal $w$-system.
\end{corollary}

Now we prove some extra properties of mincover and minpartition functions.

\begin{lemma}\label{mincover2E}
Let $f$ be a function defined on $2^E$ such that $f(\O)=0$, $\mathcal E\subseteq \mathcal A \subseteq 2^E$, and $f(Y)\geq f(Y, \mathcal A)$, for any $\O \neq Y\subseteq E$ with $Y\notin \mathcal A$. Then $f(X, 2^E)=f(X, \mathcal A)=f(X, \mathcal A_f)$, for any $\O \neq X\subseteq E$.
\end{lemma}
\begin{proof}
Let $\mathcal A_X$ be an optimal $2^E$-cover for $X$. If $\mathcal A_X\subseteq \mathcal A_f$, then we are done. Otherwise, let $\mathcal A_X\cap \mathcal A_f=\{ A_1, ..., A_p\}$ and $\mathcal A_X\backslash \mathcal A_f=\{ B_1, ..., B_q\}$. It follows that $f(X, 2^E)=\sum_{i=1}^{p} f(A_i)+\sum_{i=1}^{q} f(B_i)$. Now, for any $i=1, ..., q$, let $\mathcal B_i=\{ A_{i1}, ..., A_{ip_i}\}$ be an optimal $\mathcal A_f$-cover for $B_i$. It follows that $\mathcal A_f(X)=(\mathcal A_X\cap \mathcal A_f)\cup (\bigcup_{i=^1}^p \mathcal B_i)$ is a $\mathcal A_f$-cover of $X$. In the other hand, by relabeling the subsets of $\mathcal A_f(X)=\{ L_1, ... L_k\}$ so that we do not repeat the same subset, we have:
\\ $f(X, \mathcal A_f)\leq \sum_{j=1}^{k} f(L_j)\leq \sum_{i=1}^{p} f(A_i)+\sum_{i=1}^{p} \sum_{j=1}^{p_i} f(A_{ip_i})\leq \sum_{i=1}^{p} f(A_i)+\sum_{i=1}^{p} f(B_{i})=f(X, 2^E)$. (i) of Lemma \ref{mincover} allows us to conclude.
\end{proof}

\begin{corollary}\label{mincoverSuperL}
Let $(E, f, \mathcal A)$ be a $w$-system, $\mathcal A\subseteq \mathcal B\subseteq 2^E$, and $\O \neq X\subseteq E$. Then $f(X, \mathcal B)=f(X, \mathcal A)=f(X, \mathcal A_f)$.
\end{corollary}

\begin{corollary}\label{SuperWsystem}
Let $(E, f, \mathcal A)$ be a $w$-system, and $\mathcal A\subseteq \mathcal B\subseteq 2^E$. Then $(E, f, \mathcal B)$ is a $w$-system.
\end{corollary}

\begin{corollary}\label{fpseudosuperset}
Let $(E, f, \mathcal A)$ be a $w$-system, $A$ and $X$ two nonempty subsets. Then $f(X)\leq f(A)+|X\backslash A|$.
\end{corollary}
\begin{proof}
Direct from Lemma \ref{mincover2E} because $A\cup \mathcal E(X\backslash A)$ is an $2^E$-cover of $X$.
\end{proof}

\begin{corollary}\label{funion}
Let $(E, f, \mathcal A)$ be a $w$-system, $A$ and $B$ two nonempty subsets. Then $f(A\cup B)\leq f(A)+f(B)$ and $f(A\cup B)\leq f(A)+f(B\backslash A)$.
\end{corollary}
\begin{proof}
Direct from Lemma \ref{mincover2E} because $\{ A, B\}$ ($\{ A, B\backslash A\}$, respectively) is an $2^E$-cover of $A\cup B$.
\end{proof}

\begin{corollary}\label{fnondecreasing}
Let $(E, f, \mathcal A)$ be a $w$-system. Then
\\ (i) $f$ is nondecreasing, i.e., if $A\subseteq B\subseteq E$ than $f(A)\leq f(B)$;
\\ (ii) if $A\in \mathcal A_f$ and $A\subset B\subseteq E$ then $f(A)<f(B)$, i.e., $A$ is an $f$-flat.
\end{corollary}
\begin{proof}
Direct from Lemma \ref{mincover2E}, because $\{ B\}$ is an $2^E$-cover of $A$.
\end{proof}

\begin{corollary}\label{LfInclusion}
Let $(E, f, \mathcal L)$ be a $w$-system. Then $\mathcal L_f\subseteq \mathcal L(f)$, i.e., any subset $L\in \mathcal L_f$ is an $f$-nonseparable $f$-flat.
\end{corollary}
\begin{proof}
It suffices to prove that $L$ is $f$-nonseparable according to (ii) of Corollary \ref{fnondecreasing}. Suppose by contradiction that $L$ is $f$-separable then there exists a nonempty subset $A\subseteq L$ such that $f(L)=f(A)+f(L\backslash A)\geq f(L, 2^E)=f(L, \mathcal L)>f(L)$, a contradiction.
\end{proof}

\begin{lemma}\label{minpartition2E}
Let $(E, f, \mathcal L)$ be a $w$-system, and $\O\neq X\subseteq E$. Then at least one of the following assertions holds:
\\ (i) $f(X, \mathcal L)=\hat f(X, 2^E)$;
\\ (ii) There exists $L\in \mathcal L_f$ such that $X\subseteq L$ and $f(X, \mathcal L)=f(L)$.
\end{lemma}
\begin{proof}
Suppose that (ii) is not true. According to (iv) of Lemmas \ref{mincover} and \ref{mincover2E}, it suffices to prove that $f(X, \mathcal L)\geq \hat f(X, 2^E)$.
\\ According to Lemma \ref{mincover2E}, there exists $\mathcal A_X=\{ A_1, ..., A_p\}$ an optimal $2^E$-cover for $X$ such that $f(X, \mathcal L)=f(X, 2^E)=f(X, \mathcal A_X)$.
\\ If $p=1$ then, by choosing $A_1\supset X$ to be maximal by inclusion such that $f(X, \mathcal L)=f(A_1)$, and $A_1\notin \mathcal L_f$, we have: $f(X, \mathcal L)=\sum_{i=1}^{t} f(Z_i)$, where $t\geq 2$ and $\{ Z_1, ..., Z_t\}$ is an optimal $\mathcal L$-cover of $A_1$. So without loss of generality, we can suppose that $p\geq 2$.
\\ Since $A_i\cap X$ should be nonempty and distinct from $X$, for $i=1, ..., p$, and $f$ is nondecreasing according to Corollary \ref{fnondecreasing}, $\mathcal B_X=\{ A_1\cap X, ..., A_p\cap X\}$ is also an optimal $2^E$-cover for $X$. Now let $Y_1=A_1\cap X$, and $Y_j=(A_j\cap X)\backslash \bigcup_{i=1}^{j-1} Y_i$, $j=2, ..., p$. It is clear that $\mathcal Y_X=\{ Y_1, ... , Y_k\}$, $k\leq p$ (we keep nonempty subsets only), is an $2^{E}$-partition of $X$. Moreover, since $Y_j\subseteq A_j\cap X\subseteq A_j$ for any $j=1, ..., k$, then $f(Y_j)\leq f(A_j\cap X)\leq f(A_j)$ because $f$ is nondecreasing. Thus $\hat f(X, 2^E)\leq \sum_{i=1}^{k} f(Y_i)\leq \sum_{i=1}^{p} f(A_i)=f(X, 2^E)=f(X, \mathcal L)$.
\end{proof}

\begin{corollary}\label{nonLf}
Let $(E, f, \mathcal L)$ be a $w$-system, and $\O\neq X\subseteq E$. If $X\notin \mathcal L(f)$, then at least one of the following assertions holds:
\\ (i) There exists a partition $\{ X_1, ..., X_k\}$ of $X$ such that all $X_i$'s are $f$-nonseparable and $f(X)=\sum_{i=1}^k f(X_i)$;
\\ (ii) There exists $L\in \mathcal L(f)$ such that $X\subset L$ and $f(X)=f(L)$.
\end{corollary}

\begin{corollary}\label{Lf}
Let $(E, f, \mathcal L)$ be a $w$-system. Then $\mathcal L_f=\mathcal L(f)$, i.e., any subset $L\in \mathcal L_f$ is an $f$-nonseparable $f$-flat, and vice-versa.
\end{corollary}
\begin{proof}
  Lemma \ref{minpartition2E} implies that if $X\notin \mathcal L_f$ then it is either $f$-separable or it is not an $f$-flat, i.e., $\mathcal L(f)\subseteq \mathcal L_f$. Corollary \ref{LfInclusion} permits to conclude.
\end{proof}

\begin{corollary}\label{wIndependence}
  If $(E, f, \mathcal L)$ is a $w$-system, then $(E, \mathcal I(f))$ is an independence system.
\end{corollary}
\begin{proof}
Axiom (I1) is satisfied by axiom (L0) because $f(\O)=0=|\O|$. To prove axiom (I2), let $J\subseteq I\in \mathcal I(f)$, where $\mathcal I(f)$ is the collection of $f$-independent sets of $(E, f)$, i.e., subsets $I\subseteq E$ such that $f(I)=|I|$. According to Corollary \ref{fpseudosuperset}, $|I|=f(I)\leq f(J)+|I\backslash J|$, which means that $f(J)\geq |I|-|I\backslash J|=|J|$. (v) of Lemma \ref{mincover} allows to conclude.
\end{proof}

\begin{corollary}\label{fIndependentR3}
  Let $(E, f, \mathcal L)$ be a $w$-system, $X$ an $f$-independent set, and $Y\subseteq E$. Then: $f(X)+f(Y)\geq f(X\cap Y)+f(X\cup Y)$ ($X$ is $f$-submodular).
\end{corollary}
\begin{proof}
According to Corollary \ref{fpseudosuperset}, $f(X\cup Y)\leq f(Y)+|(X\cup Y)\backslash Y|=f(Y)+|X\backslash Y|$. This yields $f(X\cap Y)+f(X\cup Y)\leq |X\cap Y|+f(Y)+|X\backslash Y|=|X|+f(Y)=f(X)+f(Y)$.
\end{proof}

\begin{lemma}\label{laminarfSubmodular}
  Let $E$ be a finite set, $f$ a nonnegative integer function defined on $2^E$ such that $f(\O)=0$, $X$ and $Y$ two subsets. If $X\subseteq Y$ or $X\cap Y=\O$ then $(X, Y)$ is $f$-submodular, that is, $f(X\cap Y)+f(X\cup Y)\leq f(X)+f(Y)$.
\end{lemma}

\begin{corollary}\label{trivialfSubmodular}
  Let $(E, f, \mathcal L)$ be a $w$-system, and $X\in \{ \O, E\}\cup \mathcal E$. Then $X$ is $f$-submodular, that is, $f(X\cap Y)+f(X\cup Y)\leq f(X)+f(Y)$, for any $Y\subseteq E$.
\end{corollary}

\begin{lemma}\label{fCircuit}
Let $(E, f, \mathcal L)$ be a $w$-system. Then the following assertions are equivalent:
  \\ (i) $X$ is an independent set;
  \\ (ii) $X$ does not contain a circuit.
\end{lemma}
\begin{proof}
It suffice to prove that (ii) implies (i) because $(E, \mathcal I(f))$ is an independence system. We prove it by induction on $k=|X|$.
  \\ If $k=0$ then $X=\O$ and $f(X)=f(\O)=0=|X|$, which means that $X$ is an independent set.
  \\ If $k=1$ then $f(X)=1=|X|$, i.e., $X$ is an independent set, because, otherwise, $f(X)=0=|X|-1$ and the unique subset $Y\subset X$ is $Y=\O$ which is an independent set ($f(Y)=0=|Y|$), i.e., $X$ is a circuit, a contradiction with (ii).
  \\ Suppose now that $k\geq 2$. For any $e\in X$, $X_e=X\backslash \{ e\}$ does not contain a circuit. By induction on the cardinality, $X_e$ is an independent set. Hence any subset $Y\subset X$ is an independent set because there exists $e\in X$ such that $Y\subseteq X_e$ and $(E, \mathcal I(f))$ is an independence system. Suppose by contradiction that $X$ is not an independent set. It follows that $|X|-1=|X_e|=f(X_e)\leq f(X)\leq |X|-1$. In other words, $f(X)=|X|-1$, and $X$ is a circuit, a contradiction with (ii).
\end{proof}

Let $(E, f)$ be a system, and $\mathcal C\subseteq 2^E$. We introduce a new submodularity circuit axiom for $f$ and $\mathcal C$ as follows.
\\ (C4) $f(C_1)+f(C_2)\geq f(C_1\cup C_2)+f(C_1\cap C_2)$ for any $C_1, C_2\in \mathcal C$.

\begin{lemma}\label{newCircuitAxioms}
Let $(E, f, \mathcal L)$ be a $w$-system, and $\mathcal C=\mathcal C(f)$. Then circuit axioms (C1)-(C3) are equivalent to axioms (C1)-(C2) and (C4) for $f$ and $\mathcal C$.
\end{lemma}
\begin{proof}\ \\
It suffices to prove the equivalence of (C3) and (C4) under axioms (C1)-(C2).
\\ {\bf (C3) implies (C4)}
\\ Since $(C_1\cup C_2)\backslash \{ e\}$ contains a circuit, for any $e\in C_1\cap C_2$, and Lemma \ref{fCircuit}, $(C_1\cup C_2)\backslash \{ e\}$ is not an independent set and thus $f((C_1\cup C_2)\backslash \{ e\})\leq |(C_1\cup C_2)\backslash \{ e\}|-1=|C_1\cup C_2|-2$. In the other hand, $f(C_i)=|C_i|-1$, $i=1, 2$, and $f(C_1\cap C_2)=|C_1\cap C_2|$ according to Lemma \ref{fCircuit} for the latter. This means that $f(C_1)+f(C_2)-f(C_1\cap C_2)=|C_1|-1+|C_2|-1-|C_1\cap C_2|=|C_1\cup C_2|-2$.
\\ Claim: $f(C_1\cup C_2)\leq |C_1\cup C_2|-2$.
\\ Suppose that $f(C_1\cup C_2)=|C_1\cup C_2|-1$, and thus there exists $g\in C_1\cup C_2$ such that $f((C_1\cup C_2)\backslash \{ g\})=|C_1\cup C_2|-1$, i.e., $(C_1\cup C_2)\backslash \{ g\}$ is independent. In this case, if $g\notin C_1\cap C_2$, i.e., $g\in C_1\backslash C_2$ or $g\in C_1\backslash C_2$, then $(C_1\cup C_2)\backslash \{ g\}$ contains one of the two circuits $C_1$ or $C_2$, a contradiction with Lemma \ref{fCircuit}. Similarly according to axiom (C3), if $g\in C_1\cap C_2$ then $(C_1\cup C_2)\backslash \{ g\}$ contains a circuit, which contradicts that $(C_1\cup C_2)\backslash \{ g\}$ is independent.
\\ So (C3) implies (C4).
\\ {\bf (C4) implies (C3)}
Suppose, by contradiction, that there are $C_1, C_2\in \mathcal C$ and $e\in C_1\cap C_2$, such that $C_1\neq C_2$, and $(C_1\cup C_2)\backslash \{ e\}$ does not contain a circuit. It follows that $C_1\backslash \{ e\}$, $C_2\backslash \{ e\}$, and $(C_1\cup C_2)\backslash \{ e\}$ are independent sets because of Lemma \ref{fCircuit} for the latter. This means that $f(C_1\cup C_2)=|C_1\cup C_2|-1$. In the other hand, since $C_1\neq C_2$ and axiom (C2), $C_1\cap C_2\subset C_i$, $i=1, 2$. Hence $C_1\cap C_2$ is an independent set, and $f(C_1\cap C_2)=|C_1\cap C_2|$. According to (C4), $f(C_1)+f(C_2)\geq f(C_1\cup C_2)+f(C_1\cap C_2)$. This yields $|C_1|+|C_2|-2=|C_1|-1+|C_2|-1\geq |C_1\cup C_2|-1+|C_1\cap C_2|=|C_1|+|C_2|-1$, a contradiction.
\end{proof}

\begin{corollary}\label{CfMatroid}
Let $(E, f, \mathcal L)$ be a $w$-system, and $\mathcal C(f)$ satisfies axiom (C4). Then $(E, \mathcal C(f))$ is a matroid whose collection of circuits is $\mathcal C(f))$ and its rank function is $f$.
\end{corollary}
\begin{proof}
A direct consequence of Lemmas \ref{fCircuit} and \ref{newCircuitAxioms}.
\end{proof}

\begin{lemma}\label{Cfnonseparable}
  Let $(E, f, \mathcal L)$ be a $w$-system, and $C\in \mathcal C(f)$ (such that $f(C)\geq 1$). Then $C$ is $f$-nonseparable.
\end{lemma}
\begin{proof}
  Suppose by contradiction that there exists a nonempty subset $A\subset C$ such that $f(C)=f(A)+f(C\backslash A)$. By the definition of an $f$-circuit, both $A$ and $C\backslash A$ are $f$-independent, that is, $f(A)=|A|$ and $f(C\backslash A)=|C\backslash A|$. Thus $f(C)=|C|$, a contradiction.
\end{proof}

\begin{lemma}\label{fCircuitClosureW}
  Let $(E, f, \mathcal L)$ be a $w$-system, and $C\in \mathcal C(f)$ (such that $f(C)\geq 1$). Then there is an $f$-closure of $C$ that belongs to $\mathcal L_f$.
\end{lemma}
\begin{proof}
If $C$ is an $f$-flat then $C\in cl(C)$ and thus $C\in \mathcal L_f$. Otherwise, Corollary \ref{nonLf} and Lemma \ref{Cfnonseparable} imply that there exists an $f$-closure $L\supset C$ such that $f(C)=f(L)$ and $L\in \mathcal L_f$.
\end{proof}

\begin{lemma}\label{fCircuitSubmodularity}
  Let $(E, f, \mathcal L)$ be a $w_0$-system, and $C_i\in \mathcal C(f)$, $i=1, 2$. Then: $f(C_1)+f(C_2)\geq f(C_1\cup C_2)+f(C_1\cap C_2)$.
\end{lemma}
\begin{proof}
  Without loss of generality, we can suppose that $C_1\neq C_2$ and $f(C_i)\geq 1$, $i=1, 2$. Let $X_i\in cl(C_i)\cap \mathcal L_f$, $i=1, 2$, according to Lemma \ref{fCircuitClosureW}. This means that $f(C_1)+f(C_2)=f(X_1)+f(X_2)\geq f(X_1\cup X_2)+f(X_1\cap X_2)\geq f(C_1\cup C_2)+f(C_1\cap C_2)$ because $f$ is nondecreasing.
\end{proof}

\begin{theorem}\label{HLfMatroid}
  If $(E, f, \mathcal L)$ is a $w_0$-system then $(E, \mathcal C(f))$ is a matroid whose collection of circuits is $\mathcal C(f))$, its rank function is $f$, and its collection of $f$-nonseparable $f$-flats is $\mathcal L_f$.
\end{theorem}
\begin{proof}
  A direct consequence of Corollaries \ref{LfInclusion}, \ref{CfMatroid}, \ref{wIndependence}, and Lemma \ref{fCircuitSubmodularity}.
\end{proof}

\begin{theorem}\label{LfMatroid}
  If $(E, f, \mathcal L)$ is a locked $w_0$-system then $(E, \mathcal C(f))$ is a matroid whose collection of circuits is $\mathcal C(f))$, its rank function is $f$, its collection of $f$-locked subsets is $\mathcal L_{ff^*}$, and its dual is defined by $(E, f^*, \mathcal L^*)$.
\end{theorem}
\begin{proof}
  Since a locked $w_0$-system is a particular $w_0$-system then both $(E, f, \mathcal L)$ and $(E, f^*, \mathcal L^*)$ define two matroids. Actually, $(E, f, \mathcal L)$ defines a matroid $M$ and $(E, f^*, \mathcal L^*)$ defines the dual matroid $M^*$. It follows that $\mathcal L(f)=\mathcal L_{f}$ and $\mathcal L^*(f^*)=\mathcal L^*_{f^*}=(\mathcal L_{f^*})^*$. Thus $\mathcal L(ff^*)=\mathcal L_{ff^*}$.
\end{proof}

\section{Uniform matroid recognition}

We recall here some definitions and introduce some notations. Let $M$ be a matroid defined on a finite set $E$, $f$ its rank function, $M^*$ its dual, and $P\subseteq E$. We say that $P$ is a parallel closure of $M$ if $P$ is a flat of rank 1. Parallel closure are exactly $f$-nonseparable $f$-flat of rank 1. We denote by $\mathcal P(M)$ the collection of such subsets. We denote by $\mathcal L_2(M)$ the collection of $f$-nonseparable $f$-flats $L\subset E$ such that $2\leq f(L)\leq f(E)-1$, that is, $f$-nonseparable $f$-flats that are neither parallel closures nor the ground set. Finally, $\mathcal L(M)=\mathcal L(f)$. This means that $\mathcal L(M)=\mathcal P(M)\cup \mathcal L_2(M)$ if $M$ is disconnected, and $\mathcal L(M)=\{ E \}\cup \mathcal P(M)\cup \mathcal L_2(M)$ if $M$ is 2-connected.
\\ In general, we call $(E, f, \mathcal P, \mathcal L_2)$ a basic $\alpha$-quadruple if:
\\ (Q0) $f$ is defined on $\{ \O, E\}\cup \mathcal P\cup \mathcal L_2$;
\\ (Q1) $f(\O)=0$;
\\ (Q2) $0\leq f(E)\leq |E|$;
\\ (Q3) $\mathcal P$ is a partition of $E$;
\\ (Q4) $f(P)=1$ for any $P\in \mathcal P$;
\\ (Q5) $\mathcal L_2\subseteq 2^E\backslash (\{ \O, E\}\cup \mathcal P)$;
\\ (Q6) $2\leq f(L)\leq f(E)-1$ for any $L\in \mathcal L_2$;
\\ (Q7) $\alpha=2$ or $\alpha=0$ (2-connectivity or not).
\\ It is clear that the size of a basic quadruple is at most $2|E|+2|\mathcal L_2|+2$. As a consequence of this definition, we have the following corollary of Theorem \ref{HLfMatroid}.

\begin{corollary}\label{matroidInput}
 Any matroid $M$ can be described uniquely by a basic $\alpha$-quadruple satisfying axioms (L2)-(L3), where $\alpha=2$ when $M$ is 2-connected, and $\alpha=0$ otherwise.
\end{corollary}

Now we give a characterization of an uniform matroid by means of its rank and its collection of nonseparable flats.

\begin{theorem}\label{uniform}
 Let $n$ be a positive integer, $E$ a finite set of cardinality $n$, $M$ a matroid defined on $E$, and $f$ its rank function. Then $M$ is uniform if and only if one of the following properties holds:
\\ (i) $2\leq f(E)\leq |E|-2$, $|\mathcal L_2(M)|=0$, and $|\mathcal P(M)|=|E|$;
\\ (ii) $f(E)=|E|-1$, and $E\in \mathcal L(M)$;
\\ (iii) $f(E)=1$, and $E\in \mathcal L(M)$;
\\ (iv) $f(E)=|E|$;
\\ (v) $f(E)=0$.
\end{theorem}
\begin{proof}
  It is clear that we have the following equivalences:
  \\ (ii) is equivalent to $M=U_{n-1, n}$ (i.e., $E$ is a circuit (coparallel closure));
  \\ (iii) to $M=U_{1, n}$ (i.e., $E$ is a parallel closure);
  \\ (iv) to $M=U_{n, n}$ (i.e., $E$ is a basis); and
  \\ (v) to $M=U_{0, n}$ (i.e., $E$ is an union of loops (cobasis)).
\\ Suppose now that (ii)-(v) are not satisfied.
\\ If $M$ is uniform then $2\leq f(E)\leq |E|-2$. Moreover, if $L\in \mathcal L_2(M)$ then $2\leq f(L)\leq f(E)-1$, and $|L|\geq f(E)+1$. In other words, $f(L)=f(E)$, a contradiction. Thus $\mathcal L_2(M)=\O$. Furthermore, if $|\mathcal P(M)|\leq |E|-1$, then there exists a parallel closure $P$ such that $|P|\geq 2$, because $\mathcal P(M)$ is a partition of $E$. Hence $f(E)=1$, a contradiction.
\\ For the reverse, suppose that (i) holds, and let $X\subseteq E$ be a flat with $|X|\geq 2$. If $X$ is nonseparable, then $X=E$ (and $M$ is 2-connected). Otherwise, $X$ is separable into nonseparable flats of cardinality at most $|X|-1<|E|$, that is, nonseparable flats of cardinality 1. Thus $f(X)=|X|$. So any flat is either the ground set or an independent set, which means that $M$ is uniform.
\end{proof}

\begin{corollary}\label{uniformL}
 Let $n$ be a positive integer, $E$ a finite set of cardinality $n$, $M$ a matroid defined on $E$, and $f$ its rank function. Then $M$ is uniform if and only if one of the following properties holds:
\\ (i) $2\leq f(E)\leq |E|-2$, $|\mathcal L_2(M)|=0$, and $|\mathcal P(M)|=|E|$;
\\ (ii) $f(E)=|E|-1$, $|\mathcal L_2(M)|=0$, and $|\mathcal P(M)|=|E|$;
\\ (iii) $f(E)=1$, $|\mathcal L_2(M)|=0$, and $|\mathcal P(M)|=1$;
\\ (iv) $f(E)=|E|$, $|\mathcal L_2(M)|=0$, and $|\mathcal P(M)|=|E|$;
\\ (v) $f(E)=0$, $|\mathcal L_2(M)|=0$, and $|\mathcal P(M)|=0$.
\end{corollary}

For a more simplification, we call a basic $\alpha$-quadruple, a basic quadruple for short. It follows that:

\begin{corollary}\label{matroidUniform}
Let $M$ be a matroid described by its basic quadruple $(E, f, \mathcal P(M), \mathcal L_2(M))$ with $|E|\geq 2$.
\\ (i) If the size of $(E, f, \mathcal P(M), \mathcal L_2(M))$ is at least $|E|+4$ then necessarily $M$ is not uniform.
\\ (ii) Otherwise, we can decide if $M$ is uniform or not in at most $O(|E|)$ time.
\end{corollary}

\begin{corollary}\label{quadrupleUniform}
Let $(E, f, \mathcal P, \mathcal L_2)$ be a basic quadruple with $|E|\geq 2$.
\\ (i) If the size of $(E, f, \mathcal P, \mathcal L_2)$ is at least $|E|+4$ then necessarily $(E, f, \mathcal P, \mathcal L_2)$ does not define an uniform matroid.
\\ (ii) Otherwise, we can decide if $(E, f, \mathcal P, \mathcal L_2)$ defines an uniform matroid or not in at most $O(|E|)$ time.
\end{corollary}
\begin{proof}
  (i) is a consequence of Corollary \ref{matroidUniform}. For (ii), firstly we test the cases of Corollary \ref{uniformL}. If all of them are not satisfied then the basic quadruple does not define an uniform matroid. Otherwise, we can test axioms (L0)-(L1) in at most $O(|E|)$. Axioms (L2)-(L3) are necessarily satisfied. And we are done.
\end{proof}

Now we call $(E, f, \mathcal P, \mathcal L_2)$ a basic system if it satisfies axiom (Q0) only.

\begin{corollary}\label{systemUniform}
Let $(E, f, \mathcal P, \mathcal L_2)$ be a basic system with $|E|\geq 2$.
\\ (i) If the size of $(E, f, \mathcal P, \mathcal L_2)$ is at least $|E|+4$ then necessarily $(E, f, \mathcal P, \mathcal L_2)$ does not define an uniform matroid.
\\ (ii) Otherwise, we can decide if $(E, f, \mathcal P, \mathcal L_2)$ defines an uniform matroid or not in at most $O(|E|)$ time.
\end{corollary}
\begin{proof}
  (i) is a consequence of Corollary \ref{quadrupleUniform}. For (ii), firstly we test the cases of Corollary \ref{uniformL}. If all of them are not satisfied then the basic quadruple does not define an uniform matroid. Otherwise, we can test axioms (Q1)-(Q7) in at most $O(|E|)$. Actually, we have to test (Q0-(Q4) and (Q7) only, because $\mathcal L_2=\O$. If at least one of them is not satisfied, then the basic quadruple does not define an uniform matroid. Otherwise, we conclude as for the above proof of Corollary \ref{quadrupleUniform}.
\end{proof}

Note that if $|E|\geq 1$, then Corollaries \ref{matroidUniform}, \ref{quadrupleUniform}, and \ref{systemUniform} hold again.

\section{Conclusion}

Our new axioms system is worth in its own because it gives us a new point of view on how to conceive a matroid. We can have then many applications as consequences of this. As we have shown, one of these applications is recognizing uniform matroids. In the first step, our algorithms exclude any exponential input size before running the polynomial second step. Other applications will be studied in our next papers.



\vspace{11pt}
{\Large\bf No conflict of interest}

The (single) author states that there is no conflict of interest.


\end{document}